\documentclass[12pt]{article}
\usepackage[english]{babel}
\usepackage{latexsym}
\usepackage{amssymb,amsthm,amsmath}
\usepackage{graphicx}
\title{\bf A computer based classification of caps in $PG(5,2)$}
\author{Daniele Bartoli, Stefano Marcugini, Fernanda Pambianco}

{\theoremstyle{definition}
\newtheorem*{definition*}{Definition}
}

\newtheorem*{proposition*}{Proposition}
\newtheorem*{corollary*}{Corollary}
\newtheorem*{lemma*}{Lemma}
\date{}
\begin{document}

\maketitle \vspace*{-15mm} \noindent

 \vspace*{2cm}
\begin{abstract}
\noindent In this paper we present the complete classification of
caps in $PG(5,2)$. These results have been obtained using a computer
based exhaustive search that exploits projective equivalence.
\end{abstract}

\section{Introduction}
In the projective space $PG(r,q)$ over the Galois Field $GF(q)$, a $n$-cap is a set of $n$ points
 no 3 of which are collinear. A $n$-cap is called \emph{complete} if it is not contained in a
 $(n+1)$-cap. For  a detailed description of the most important properties of these geometric
 structures, we refer the reader to \cite{Hir}.
 In the last decades the problem of determining the spectrum of the
 sizes of complete caps has been the subject of a lot of researches.
 For a survey see \cite{18secondo}.\\
In this work we search for the classification of complete and
incomplete caps in $PG(5,2)$, using an exhaustive search algorithm. In Section \ref{alg}
the algorithm utilized is illustrated;
 in Section \ref{risultati} the complete list of non equivalent complete and incomplete
caps is presented.

\section{The searching algorithm }\label{alg}
In this section the algorithm utilized is presented. Our goal is to obtain the classification of complete and incomplete caps in $PG(5,2)$. It is not restrictive to suppose that a cap in $PG(5,2)$ contains this six points:
\begin{displaymath}
\mathcal{R}=\{(1:0:0:0:0:0);(0:1:0:0:0:0);(0:0:1:0:0:0);
\end{displaymath}
\begin{displaymath}
(0:0:0:1:0:0);(0:0:0:0:1:0);(0:0:0:0:0:1)\}.
\end{displaymath}
Then we define the set \emph{Cand} of all the points lying no 2-secant of $\mathcal{R}$. We introduce in \emph{Cand} the following equivalence relationship:
\begin{displaymath}
P \sim Q \iff \mathcal{C} \cup \{P\} \cong \mathcal{C} \cup \{Q\},
\end{displaymath}

\noindent where $\cong$ means that the two sets are projectively equivalent. This relationship spreads the candidates in equivalent classes $\mathcal{C}_{1},\ldots,\mathcal{C}_{k}$.\\
The choice of the next point to add to the building cap can be made only among the representatives of the equivalent classes, in fact two caps one containing $\mathcal{C} \cup \{P\}$ and the other one $\mathcal{C} \cup \{Q\}$, with $P$ and $Q$ in $\mathcal{C}_{\overline{i}}$, are equivalent by definition of orbit.\\
Suppose now that we have construct all the caps containing $\mathcal{C} \cup \{P_{i}\}$, with $i \leq \overline{i}$. Considering the caps containing $\mathcal{C} \cup \{P_{j}\}$ with $\overline{i}<j$, all the points of the classes $\mathcal{C}_{k}$ with $k<\overline{i}$ can be avoided. In fact a cap containing $\mathcal{C} \cup \{P_{j}\} \cup \{\overline{P}_{k}\}$, with $\overline{P}_{k} \in \mathcal{C}_{k}$ and $k<\overline{i}$, is projectively equivalent to a cap containing $\mathcal{C} \cup \{P_{k}\} \cup \{P_{j}\}$, already studied.\\
When we add a new point to the cap, we can divide all the remaining candidates in equivalence classes, as above. Two points  $P$ and $Q$ are in relationship with the $j$-th class $\mathcal{C}_{j}$, i.e. $P \sim_{j} Q$, if $\mathcal{C} \cup \{P_{j}\} \cup \{P\}$ and $\mathcal{C} \cup \{P_{j}\} \cup \{Q\}$ are projectively equivalent. \\
At the $m$-th step of the extension process if the cap $\mathcal{C}\cup \{P_{i_{m}}\} \cup \ldots \cup \{P^{i_{1}\ldots i_{m-1}}_{i_{m}}\} \cup\{P\}$  is projectively equivalent to the cap $\mathcal{C}\cup \{P_{i_{m}}\} \cup \ldots \cup \{P^{i_{1}\ldots i_{m-1}}_{i_{m}}\} \cup\{Q\}$ with $P^{i_{1}\ldots i_{r}}_{s} \in \mathcal{C}^{i_{1}\ldots i_{r}}_{s} $, then $P$ and $Q$ are in relationship ($P \sim_{i_{1}\ldots i_{m}} Q$) and they belong to the same class $\mathcal{C}^{i_{1}\ldots i_{m}}_{m+1} $.\\
Iterating the process we can build a tree similar to the following:
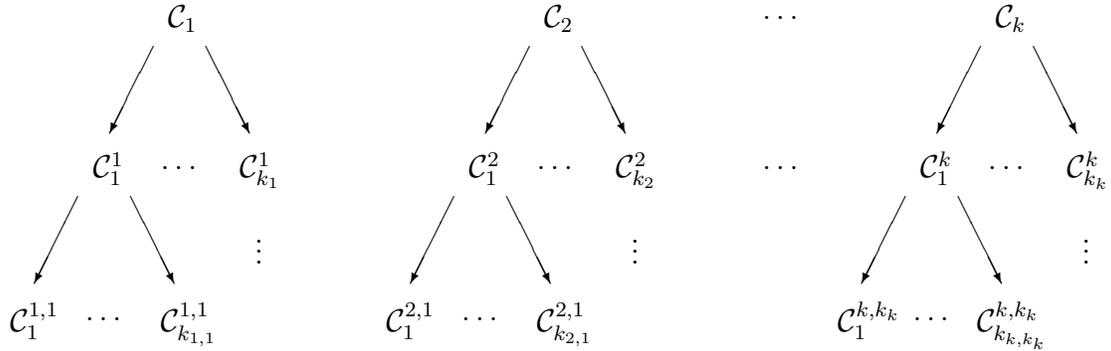
\begin{center}

\setlength{\unitlength}{1cm}
\begin{picture}(15,6)(0,0)

%Primo blocco

\put(2.35,4.9){$\mathcal{C}_{1}$}

\put(2.15,4.65){\vector(-1,-2){0.58}}
\put(2.85,4.65){\vector(1,-2){0.58}}

\put(1.35,2.9){$\mathcal{C}_{1}^{1}$}
\put(2.25,3){\ldots}

\put(3.3,2.9){$\mathcal{C}_{k_{1}}^{1}$}

\put(1.15,2.65){\vector(-1,-2){0.58}}
\put(1.85,2.65){\vector(1,-2){0.58}}

\put(3.5,1.75){\vdots}

\put(0.25,0.85){$\mathcal{C}_{1}^{1,1}$}
\put(1.25,1){\ldots}

\put(2.25,0.85){$\mathcal{C}_{k_{1,1}}^{1,1}$}

\put(10.25,5){\ldots}
\put(10.25,3){\ldots}

%Secondo blocco

\put(7.35,4.9){$\mathcal{C}_{2}$}

\put(7.15,4.65){\vector(-1,-2){0.58}}
\put(7.85,4.65){\vector(1,-2){0.58}}

\put(6.35,2.9){$\mathcal{C}_{1}^{2}$}
\put(7.25,3){\ldots}

\put(8.30,2.9){$\mathcal{C}_{k_{2}}^{2}$}

\put(6.15,2.65){\vector(-1,-2){0.58}}
\put(6.85,2.65){\vector(1,-2){0.58}}

\put(8.5,1.75){\vdots}

\put(5.25,0.85){$\mathcal{C}_{1}^{2,1}$}
\put(6.25,1){\ldots}

\put(7.25,0.85){$\mathcal{C}_{k_{2,1}}^{2,1}$}

%Terzo blocco

\put(13.35,4.9){$\mathcal{C}_{k}$}

\put(13.15,4.65){\vector(-1,-2){0.58}}
\put(13.85,4.65){\vector(1,-2){0.58}}

\put(12.35,2.9){$\mathcal{C}_{1}^{k}$}
\put(13.25,3){\ldots}

\put(14.3,2.9){$\mathcal{C}_{k_{k}}^{k}$}

\put(12.15,2.65){\vector(-1,-2){0.58}}
\put(12.85,2.65){\vector(1,-2){0.58}}

\put(14.5,1.75){\vdots}

\put(11.25,0.85){$\mathcal{C}_{1}^{k,k_{k}}$}
\put(12.25,1){\ldots}

\put(13.1,0.85){$\mathcal{C}_{k_{k,k_{k}}}^{k,k_{k}}$}

\end{picture}
\end{center}
The tree is important to restrict the number of candidates in the extension process. Suppose that we have generated a $n$-cap containing the cap $\mathcal{C}\cup \{P_{i_{1}}\} \cup \{P_{i_{2}}^{i_{1}}\} \cup \ldots \cup \{P^{i_{1}\ldots i_{m-1}}_{i_{m}}\} \cup\{P\}$, after having generated $n$-caps containing $\mathcal{C}\cup \{P_{j}\}$ with $j<i_{1}$, $\mathcal{C}\cup \{P_{i_{1}}\}\cup \{P^{i_{1}}_{j}\}$ with $j<i_{2}$,\ldots, $\mathcal{C}\cup \{P_{i_{1}}\}\cup \{P^{i_{1}}_{i_{2}}\}\cup \ldots \cup \{P^{i_{1}\ldots i_{m-1}}_{j}\}$ with $j<i_{m}$, with $P^{i_{1}\ldots i_{r}}_{s} \in \mathcal{C}^{i_{1}\ldots i_{r}}_{s}$. Then the points belonging to $\mathcal{C}_{1} \cup\ldots \cup \mathcal{C}_{i_{1}-1} \cup\mathcal{C}^{i_{1}}_{1} \cup\ldots \cup \mathcal{C}^{i_{1}}_{i_{2}-1} \cup \ldots \cup \mathcal{C}^{i_{1}\ldots i_{m-1}}_{1} \cup\ldots \cup \mathcal{C}^{i_{1}\ldots i_{m-1}}_{i_{m}-1}$ can be avoided, because a cap containing one of them is equivalent to one already found. For example a $n$-cap containing $\mathcal{C}\cup \{P_{i_{1}}\} \cup \ldots \cup \{P^{i_{1}\ldots i_{m-1}}_{i_{m}}\} \cup\{P\}\cup \{Q\}$ with $Q \in \mathcal{C}_{h}$ for some $h<i_{1}$ is equivalent to a $n$-cap containing $\mathcal{C} \cup \{P_{h}\}$, which is already found. \\

\section{Results}\label{risultati}
In this Section all non equivalent caps, complete and incomplete, in
$PG(5,2)$ are presented. 
\subsection{Non-equivalent caps $\mathcal{K}$ in $PG(5,2)$}
This table shows the number and the type of the non equivalent
examples of all the caps.
\begin{table}[h]
\caption{Number and type of non equivalent examples}
\begin{center}
\begin{tabular}{|c|c|c||c|c|c|}
\hline
$|\mathcal{K}|$&\# COMPLETE&\# INCOMPLETE&$|\mathcal{K}|$&\# COMPLETE&\# INCOMPLETE\\
&CAPS&CAPS&&CAPS&CAPS\\
\hline
7&0&4&20&1&23\\
8&0&7&21&0&16\\
9&0&12&22&0&15\\
10&0&24&23&0&9\\
11&0&34&24&0&8\\
12&0&43&25&0&5\\
13&1&46&26&0&4\\
14&0&49&27&0&2\\
15&0&44&28&0&2\\
16&0&48&29&0&1\\
17&5&35&30&0&1\\
18&1&32&31&0&1\\
19&0&25&32&1&0\\
\hline
\end{tabular}
\end{center}
\end{table}

\subsection{Description of the caps}
In this section we describe each cap of size $k+1$, with $7 \geq k \geq 31$, as union of a cap of size $k$ and a point in $PG(5,2)$. 
We start from these non equivalent examples of $7$-incomplete caps.
\begin{minipage}[t]{7 cm}
\begin{center}
CAP 1
\end{center}
\begin{center}
\begin{tabular}{ccccccc}
 1& 1& 0& 0& 0& 0& 0\\
 0& 1& 1& 0& 0& 0& 0\\
 0& 1& 0& 1& 0& 0& 0\\
 0& 1& 0& 0& 1& 0& 0\\
 0& 1& 0& 0& 0& 1& 0\\
 0& 1& 0& 0& 0& 0& 1\\
\end{tabular}
\end{center}
 \end{minipage}
 \begin{minipage}[t]{2,5 cm}
 $\textrm{}$\\
 \end{minipage}
\begin{minipage}[t]{7 cm}
\begin{center}
CAP 2
\end{center}
\begin{center}
\begin{tabular}{ccccccc}
 1& 1& 0& 0& 0& 0& 0\\
 0& 0& 1& 0& 0& 0& 0\\
 0& 0& 0& 1& 0& 0& 0\\
 0& 0& 0& 0& 1& 0& 0\\
 1& 0& 0& 0& 0& 1& 0\\
 1& 0& 0& 0& 0& 0& 1\\
\end{tabular}
\end{center}
 \end{minipage}
$\textrm{}$\\$\textrm{}$\\
\begin{minipage}[t]{7 cm}
\begin{center}
CAP 3
\end{center}
\begin{center}
\begin{tabular}{ccccccc}
 1& 0& 0& 0& 0& 0& 1\\
 0& 1& 0& 0& 0& 0& 1\\
 0& 0& 1& 0& 0& 0& 0\\
 0& 0& 0& 1& 0& 0& 1\\
 0& 0& 0& 0& 1& 0& 1\\
 0& 0& 0& 0& 0& 1& 0\\
\end{tabular}
\end{center}
 \end{minipage}
 \begin{minipage}[t]{2,5 cm}
 $\textrm{}$\\
 \end{minipage}
\begin{minipage}[t]{7 cm}
\begin{center}
CAP 4
\end{center}
\begin{center}
\begin{tabular}{ccccccc}
 1& 0& 0& 0& 0& 0& 1\\
 0& 1& 0& 0& 0& 0& 1\\
 0& 0& 1& 0& 0& 0& 1\\
 0& 0& 0& 1& 0& 0& 0\\
 0& 0& 0& 0& 1& 0& 1\\
 0& 0& 0& 0& 0& 1& 1\\
\end{tabular}
\end{center}
\end{minipage}

$\textrm{ }$\\
$\textrm{ }$\\
We will call these caps $\mathcal{C}_{7}^{1}, \mathcal{C}_{7}^{2}, \mathcal{C}_{7}^{3}, \mathcal{C}_{7}^{4}$ respectively. They have stabilizer ($\mathcal{O}$) of size
 
\begin{table}[h]
\caption{Cap $\mathcal{K}$ with $|\mathcal{K}|=7$}
\begin{center}
\begin{tabular}{|c|c|}
\hline
$\mathcal{K}$&STABILIZER\\
\hline
\hline
$\mathcal{C}_{7}^{1}$&$|\mathcal{O}|=5040$\\
$\mathcal{C}_{7}^{2}$&$|\mathcal{O}|=144$\\
$\mathcal{C}_{7}^{3}$&$|\mathcal{O}|=240$\\
$\mathcal{C}_{7}^{4}$&$|\mathcal{O}|=720$\\
\hline
\end{tabular}
\end{center}
\end{table}
In the following tables the symbol $\mathcal{C}_{i}^{j}$ indicates the $j$-th cap of size $i$ and  $\overline{\mathcal{C}}_{i}^{j}$ means that the cap is complete. If it is possible to write a cap in two different ways, we choose that one with the lowest value of $j$.
%\begin{table}[h]
%\caption{Cap $\mathcal{K}$ with $|\mathcal{K}|=8$}
\begin{center}
\begin{tabular}{|c|c|c||c|c|c|}
\hline
$\mathcal{K}$&CORRESPONDS TO&STABILIZER&$\mathcal{K}$&CORRESPONDS TO&STABILIZER\\
\hline
\hline
$\mathcal{C}_{8}^{1}$&$\mathcal{C}_{7}^{2}\cup\{(1,1,0,1,1,0)\}$&$|\mathcal{O}|=48$&
$\mathcal{C}_{8}^{2}$&$\mathcal{C}_{7}^{3}\cup\{(0,1,1,0,1,1)\}$&$|\mathcal{O}|=144$\\
$\mathcal{C}_{8}^{3}$&$\mathcal{C}_{7}^{2}\cup\{(1,0,0,1,0,1)\}$&$|\mathcal{O}|=96$&
$\mathcal{C}_{8}^{4}$&$\mathcal{C}_{7}^{2}\cup\{(1,1,1,1,1,0)\}$&$|\mathcal{O}|=192$\\
$\mathcal{C}_{8}^{5}$&$\mathcal{C}_{7}^{2}\cup\{(0,1,1,1,0,0)\}$&$|\mathcal{O}|=1152$&
$\mathcal{C}_{8}^{6}$&$\mathcal{C}_{7}^{2}\cup\{(1,1,1,0,1,1)\}$&$|\mathcal{O}|=72$\\
$\mathcal{C}_{8}^{7}$&$\mathcal{C}_{7}^{1}\cup\{(1,1,0,1,1,0)\}$&$|\mathcal{O}|=144$&&&\\
\hline
\end{tabular}
\end{center}
%\end{table}
%\begin{table}[h]
%\caption{Cap $\mathcal{K}$ with $|\mathcal{K}|=9$}
\begin{center}
\begin{tabular}{|c|c|c||c|c|c|}
\hline
$\mathcal{K}$&CORRESPONDS TO&STABILIZER&$\mathcal{K}$&CORRESPONDS TO&STABILIZER\\
\hline
\hline
$\mathcal{C}_{9}^{1}$&$\mathcal{C}_{8}^{1}\cup\{(0,0,1,0,1,1)\}$&$|\mathcal{O}|=48$&
$\mathcal{C}_{9}^{2}$&$\mathcal{C}_{8}^{1}\cup\{(1,0,0,1,0,1)\}$&$|\mathcal{O}|=48$\\
$\mathcal{C}_{9}^{3}$&$\mathcal{C}_{8}^{3}\cup\{(1,1,0,0,1,0)\}$&$|\mathcal{O}|=32$&
$\mathcal{C}_{9}^{4}$&$\mathcal{C}_{8}^{4}\cup\{(1,1,1,0,1,1)\}$&$|\mathcal{O}|=48$\\
$\mathcal{C}_{9}^{5}$&$\mathcal{C}_{8}^{3}\cup\{(1,1,1,0,1,1)\}$&$|\mathcal{O}|=48$&
$\mathcal{C}_{9}^{6}$&$\mathcal{C}_{8}^{3}\cup\{(1,0,1,0,0,1)\}$&$|\mathcal{O}|=384$\\
$\mathcal{C}_{9}^{7}$&$\mathcal{C}_{8}^{1}\cup\{(1,1,1,0,1,0)\}$&$|\mathcal{O}|=128$&
$\mathcal{C}_{9}^{8}$&$\mathcal{C}_{8}^{1}\cup\{(1,0,1,0,1,0)\}$&$|\mathcal{O}|=288$\\
$\mathcal{C}_{9}^{9}$&$\mathcal{C}_{8}^{1}\cup\{(0,1,1,0,1,1)\}$&$|\mathcal{O}|=32$&
$\mathcal{C}_{9}^{10}$&$\mathcal{C}_{8}^{3}\cup\{(0,0,0,1,1,1)\}$&$|\mathcal{O}|=336$\\
$\mathcal{C}_{9}^{11}$&$\mathcal{C}_{8}^{3}\cup\{(1,1,1,1,1,0)\}$&$|\mathcal{O}|=144$&
$\mathcal{C}_{9}^{12}$&$\mathcal{C}_{8}^{1}\cup\{(1,1,1,0,1,1)\}$&$\mathcal{D}_{4}\times \mathcal{C}_{2}$\\
\hline
\end{tabular}
\end{center}
%\end{table}
%\begin{table}[h]
%\caption{Cap $\mathcal{K}$ with $|\mathcal{K}|=10$}
\begin{center}
\begin{tabular}{|c|c|c||c|c|c|}
\hline
$\mathcal{K}$&CORRESPONDS TO&STABILIZER&$\mathcal{K}$&CORRESPONDS TO&STABILIZER\\
\hline
\hline
$\mathcal{C}_{10}^{1}$&$\mathcal{C}_{9}^{2}\cup\{(0,0,0,1,1,1)\}$&$|\mathcal{O}|=336$&
$\mathcal{C}_{10}^{2}$&$\mathcal{C}_{9}^{3}\cup\{(0,0,1,0,1,1)\}$&$|\mathcal{O}|=192$\\
$\mathcal{C}_{10}^{3}$&$\mathcal{C}_{9}^{2}\cup\{(0,1,1,1,1,1)\}$&$|\mathcal{O}|=48$&
$\mathcal{C}_{10}^{4}$&$\mathcal{C}_{9}^{3}\cup\{(0,1,1,1,0,0)\}$&$\mathcal{D}_{4}\times \mathcal{C}_{2}$\\
$\mathcal{C}_{10}^{5}$&$\mathcal{C}_{9}^{2}\cup\{(1,0,1,1,1,0)\}$&$|\mathcal{O}|=192$&
$\mathcal{C}_{10}^{6}$&$\mathcal{C}_{9}^{1}\cup\{(1,1,1,0,1,1)\}$&$\mathcal{D}_{4}\times \mathcal{C}_{2}$\\
$\mathcal{C}_{10}^{7}$&$\mathcal{C}_{9}^{6}\cup\{(1,1,0,0,0,1)\}$&$|\mathcal{O}|=3840$&
$\mathcal{C}_{10}^{8}$&$\mathcal{C}_{9}^{3}\cup\{(0,1,0,1,1,0)\}$&$|\mathcal{O}|=72$\\
$\mathcal{C}_{10}^{9}$&$\mathcal{C}_{9}^{10}\cup\{(1,0,0,1,1,0)\}$&$|\mathcal{O}|=2688$&
$\mathcal{C}_{10}^{10}$&$\mathcal{C}_{9}^{7}\cup\{(0,1,1,1,1,1)\}$&$|\mathcal{O}|=64$\\
$\mathcal{C}_{10}^{11}$&$\mathcal{C}_{9}^{1}\cup\{(1,0,1,0,0,1)\}$&$|\mathcal{O}|=144$&
$\mathcal{C}_{10}^{12}$&$\mathcal{C}_{9}^{3}\cup\{(0,0,0,1,1,1)\}$&$|\mathcal{O}|=48$\\
$\mathcal{C}_{10}^{13}$&$\mathcal{C}_{9}^{3}\cup\{(1,1,1,0,1,1)\}$&$\mathcal{D}_{4}\times \mathcal{C}_{2}$&
$\mathcal{C}_{10}^{14}$&$\mathcal{C}_{9}^{4}\cup\{(1,1,0,1,0,1)\}$&$\mathcal{D}_{6}$\\
\hline
\end{tabular}
\end{center}
%\end{table}

%\newpage

%\begin{table}[h]
%\caption{Cap $\mathcal{K}$ with $|\mathcal{K}|=11$}
\begin{center}
\begin{tabular}{|c|c|c||c|c|c|}
\hline
$\mathcal{K}$&CORRESPONDS TO&STABILIZER&$\mathcal{K}$&CORRESPONDS TO&STABILIZER\\
\hline
\hline
$\mathcal{C}_{11}^{1}$&$\mathcal{C}_{10}^{6}\cup\{(0,1,1,0,1,0)\}$&$|\mathcal{O}|=144$&
$\mathcal{C}_{11}^{2}$&$\mathcal{C}_{10}^{6}\cup\{(0,1,1,0,0,1)\}$&$\mathcal{D}_{6} \times \mathcal{C}_{2}$\\
$\mathcal{C}_{11}^{3}$&$\mathcal{C}_{10}^{9}\cup\{(0,1,1,0,1,1)\}$&$|\mathcal{O}|=1152$&
$\mathcal{C}_{11}^{4}$&$\mathcal{C}_{10}^{4}\cup\{(1,0,1,0,0,1)\}$&$|\mathcal{O}|=32$\\
$\mathcal{C}_{11}^{5}$&$\mathcal{C}_{10}^{2}\cup\{(0,1,1,1,0,1)\}$&$|\mathcal{O}|=48$&
$\mathcal{C}_{11}^{6}$&$\mathcal{C}_{10}^{3}\cup\{(1,1,1,0,0,1)\}$&$|\mathcal{O}|=96$\\
$\mathcal{C}_{11}^{7}$&$\mathcal{C}_{10}^{9}\cup\{(1,1,1,0,1,1)\}$&$|\mathcal{O}|=1008$&
$\mathcal{C}_{11}^{8}$&$\mathcal{C}_{10}^{10}\cup\{(0,0,1,1,1,0)\}$&$\mathcal{D}_{4} \times \mathcal{C}_{2}$\\
$\mathcal{C}_{11}^{9}$&$\mathcal{C}_{10}^{9}\cup\{(1,1,0,0,1,0)\}$&$|\mathcal{O}|=1384$&
$\mathcal{C}_{11}^{10}$&$\mathcal{C}_{10}^{14}\cup\{(1,0,1,1,0,0)\}$&$|\mathcal{O}|=120$\\
$\mathcal{C}_{11}^{11}$&$\mathcal{C}_{10}^{3}\cup\{(1,1,0,0,0,1)\}$&$|\mathcal{O}|=384$&
$\mathcal{C}_{11}^{12}$&$\mathcal{C}_{10}^{5}\cup\{(0,1,1,0,1,1)\}$&$|\mathcal{O}|=32$\\
$\mathcal{C}_{11}^{13}$&$\mathcal{C}_{10}^{8}\cup\{(0,0,0,1,1,1)\}$&$|\mathcal{O}|=48$&
$\mathcal{C}_{11}^{14}$&$\mathcal{C}_{10}^{3}\cup\{(1,1,1,0,1,1)\}$&$\mathcal{D}_{4}$\\
\hline
\end{tabular}
\end{center}
%\end{table}

\newpage

%\begin{table}[h]
%\caption{Cap $\mathcal{K}$ with $|\mathcal{K}|=12$}
\begin{center}
\begin{tabular}{|c|c|c||c|c|c|}
\hline
$\mathcal{K}$&CORRESPONDS TO&STABILIZER&$\mathcal{K}$&CORRESPONDS TO&STABILIZER\\
\hline
\hline
$\mathcal{C}_{12}^{1}$&$\mathcal{C}_{11}^{6}\cup\{(0,0,1,1,1,0)\}$&$|\mathcal{O}|=32$&
$\mathcal{C}_{12}^{2}$&$\mathcal{C}_{11}^{13}\cup\{(1,1,1,0,1,1)\}$&$\mathcal{D}_{4}$\\
$\mathcal{C}_{12}^{3}$&$\mathcal{C}_{11}^{2}\cup\{(1,0,1,0,0,1)\}$&$|\mathcal{O}|=48$&
$\mathcal{C}_{12}^{4}$&$\mathcal{C}_{11}^{10}\cup\{(1,1,0,0,1,0)\}$&$\mathcal{D}_{5}$\\
$\mathcal{C}_{12}^{5}$&$\mathcal{C}_{11}^{4}\cup\{(0,1,1,1,1,1)\}$&$|\mathcal{O}|=128$&
$\mathcal{C}_{12}^{6}$&$\mathcal{C}_{11}^{6}\cup\{(1,0,1,1,1,0)\}$&$|\mathcal{O}|=1152$\\
$\mathcal{C}_{12}^{7}$&$\mathcal{C}_{11}^{34}\cup\{(0,1,1,1,1,1)\}$&$|\mathcal{O}|=23040$&
$\mathcal{C}_{12}^{8}$&$\mathcal{C}_{11}^{5}\cup\{(1,1,1,1,1,1)\}$&$|\mathcal{O}|=64$\\
$\mathcal{C}_{12}^{9}$&$\mathcal{C}_{11}^{27}\cup\{(0,1,1,1,0,0)\}$&$\mathcal{D}_{6}$&
$\mathcal{C}_{12}^{10}$&$\mathcal{C}_{11}^{4}\cup\{(1,1,1,1,1,0)\}$&$|\mathcal{O}|=384$\\
$\mathcal{C}_{12}^{11}$&$\mathcal{C}_{11}^{1}\cup\{(0,1,1,0,0,1)\}$&$\mathcal{D}_{6} \times \mathcal{C}_{2}$&
$\mathcal{C}_{12}^{12}$&$\mathcal{C}_{11}^{7}\cup\{(1,1,0,0,1,0)\}$&$|\mathcal{O}|=96$\\
$\mathcal{C}_{12}^{13}$&$\mathcal{C}_{11}^{3}\cup\{(1,1,0,0,1,0)\}$&$|\mathcal{O}|=128$&
$\mathcal{C}_{12}^{14}$&$\mathcal{C}_{11}^{13}\cup\{(1,1,1,0,0,0)\}$&$\mathcal{D}_{4} \times \mathcal{C}_{2}\times \mathcal{C}_{2}$\\
\hline
\end{tabular}
\end{center}
%\end{table}

\newpage

%\begin{table}[h]
%\caption{Cap $\mathcal{K}$ with $|\mathcal{K}|=13$}
\begin{center}
\begin{tabular}{|c|c|c||c|c|c|}
\hline
$\mathcal{K}$&CORRESPONDS TO&STABILIZER&$\mathcal{K}$&CORRESPONDS TO&STABILIZER\\
\hline
\hline
$\overline{\mathcal{C}}_{13}^{1}$&$\mathcal{C}_{12}^{6}\cup\{(0,0,0,1,1,1)\}$&$|\mathcal{O}|=1152$&
$\mathcal{C}_{13}^{2}$&$\mathcal{C}_{12}^{2}\cup\{(1,0,1,0,0,1)\}$&$\mathcal{D}_{6}$\\
$\mathcal{C}_{13}^{3}$&$\mathcal{C}_{12}^{12}\cup\{(1,1,1,0,0,0)\}$&$|\mathcal{O}|=32$&
$\mathcal{C}_{13}^{4}$&$\mathcal{C}_{12}^{28}\cup\{(1,1,1,1,0,0)\}$&$|\mathcal{O}|=48$\\
$\mathcal{C}_{13}^{5}$&$\mathcal{C}_{12}^{2}\cup\{(1,1,1,0,0,0)\}$&$\mathcal{C}_{2} \times \mathcal{C}_{2}$&
$\mathcal{C}_{13}^{6}$&$\mathcal{C}_{12}^{11}\cup\{(0,1,1,1,0,0)\}$&$\mathcal{C}_{2} \times \mathcal{C}_{2} \times \mathcal{C}_{2}$\\
$\mathcal{C}_{13}^{7}$&$\mathcal{C}_{12}^{11}\cup\{(1,1,1,0,0,0)\}$&$|\mathcal{O}|=96$&
$\mathcal{C}_{13}^{8}$&$\mathcal{C}_{12}^{14}\cup\{(0,1,1,1,0,0)\}$&$\mathcal{D}_{4} \times \mathcal{C}_{2}$\\
$\mathcal{C}_{13}^{9}$&$\mathcal{C}_{12}^{1}\cup\{(0,0,0,1,1,1)\}$&$|\mathcal{O}|=32$&
$\mathcal{C}_{13}^{10}$&$\mathcal{C}_{12}^{2}\cup\{(0,1,0,1,0,1)\}$&$|\mathcal{O}|=48$\\
$\mathcal{C}_{13}^{11}$&$\mathcal{C}_{12}^{23}\cup\{(1,0,0,1,1,0)\}$&$|\mathcal{O}|=192$&
$\mathcal{C}_{13}^{12}$&$\mathcal{C}_{12}^{27}\cup\{(0,1,1,1,0,1)\}$&$\mathcal{D}_{6}$\\
$\mathcal{C}_{13}^{13}$&$\mathcal{C}_{12}^{3}\cup\{(1,0,1,0,1,0)\}$&$|\mathcal{O}|=384$&
$\mathcal{C}_{13}^{14}$&$\mathcal{C}_{12}^{15}\cup\{(1,1,0,0,0,1)\}$&$|\mathcal{O}|=64$\\
\hline
\end{tabular}
\end{center}
%\end{table}

\newpage

%\begin{table}[h]
%\caption{Cap $\mathcal{K}$ with $|\mathcal{K}|=14$}
\begin{center}
\begin{tabular}{|c|c|c||c|c|c|}
\hline
$\mathcal{K}$&CORRESPONDS TO&STABILIZER&$\mathcal{K}$&CORRESPONDS TO&STABILIZER\\
\hline
\hline
$\mathcal{C}_{14}^{1}$&$\mathcal{C}_{13}^{16}\cup\{(1,1,1,0,1,1)\}$&$|\mathcal{O}|=48$&
$\mathcal{C}_{14}^{2}$&$\mathcal{C}_{13}^{3}\cup\{(0,1,1,0,0,1)\}$&$\mathcal{D}_{4} \times \mathcal{C}_{2}$\\
$\mathcal{C}_{14}^{3}$&$\mathcal{C}_{13}^{22}\cup\{(1,0,0,1,1,0)\}$&$|\mathcal{O}|=128$&
$\mathcal{C}_{14}^{4}$&$\mathcal{C}_{13}^{2}\cup\{(0,1,1,0,1,0)\}$&$|\mathcal{O}|=48$\\
$\mathcal{C}_{14}^{5}$&$\mathcal{C}_{13}^{7}\cup\{(1,0,1,0,0,1)\}$&$|\mathcal{O}|=96$&
$\mathcal{C}_{14}^{6}$&$\mathcal{C}_{13}^{6}\cup\{(1,1,0,1,0,1)\}$&$\mathcal{D}_{4} \times \mathcal{C}_{2}$\\
$\mathcal{C}_{14}^{7}$&$\mathcal{C}_{13}^{5}\cup\{(0,1,1,1,0,0)\}$&$\mathcal{C}_{2} \times \mathcal{C}_{2}$&
$\mathcal{C}_{14}^{8}$&$\mathcal{C}_{13}^{6}\cup\{(0,0,0,1,1,1)\}$&$|\mathcal{O}|=36$\\
$\mathcal{C}_{14}^{9}$&$\mathcal{C}_{13}^{4}\cup\{(0,0,1,1,1,1)\}$&$|\mathcal{O}|=192$&
$\mathcal{C}_{14}^{10}$&$\mathcal{C}_{13}^{9}\cup\{(1,0,1,1,0,0)\}$&$\mathcal{D}_{6} \times \mathcal{C}_{2}$\\
$\mathcal{C}_{14}^{11}$&$\mathcal{C}_{13}^{3}\cup\{(0,0,1,1,1,0)\}$&$|\mathcal{O}|=48$&
$\mathcal{C}_{14}^{12}$&$\mathcal{C}_{13}^{6}\cup\{(1,1,1,0,0,0)\}$&$|\mathcal{O}|=32$\\
$\mathcal{C}_{14}^{13}$&$\mathcal{C}_{13}^{8}\cup\{(1,0,1,1,0,0)\}$&$|\mathcal{O}|=128$&
$\mathcal{C}_{14}^{14}$&$\mathcal{C}_{13}^{5}\cup\{(0,0,1,0,1,1)\}$&$\mathcal{D}_{4} \times \mathcal{C}_{2}$\\
\hline
\end{tabular}
\end{center}
%\end{table}

\newpage

%\begin{table}[h]
%\caption{Cap $\mathcal{K}$ with $|\mathcal{K}|=15$}
\begin{center}
\begin{tabular}{|c|c|c||c|c|c|}
\hline
$\mathcal{K}$&CORRESPONDS TO&STABILIZER&$\mathcal{K}$&CORRESPONDS TO&STABILIZER\\
\hline
\hline
$\mathcal{C}_{15}^{1}$&$\mathcal{C}_{14}^{19}\cup\{(0,1,1,0,0,1)\}$&$\mathcal{D}_{4}$&
$\mathcal{C}_{15}^{2}$&$\mathcal{C}_{14}^{4}\cup\{(0,1,1,0,0,1)\}$&$\mathcal{C}_{2} \times \mathcal{C}_{2}$\\
$\mathcal{C}_{15}^{3}$&$\mathcal{C}_{14}^{7}\cup\{(0,0,1,1,1,0)\}$&$\mathcal{S}_{3}$&
$\mathcal{C}_{15}^{4}$&$\mathcal{C}_{14}^{7}\cup\{(1,0,1,0,0,1)\}$&$\mathcal{D}_{4} \times \mathcal{C}_{2}$\\
$\mathcal{C}_{15}^{5}$&$\mathcal{C}_{14}^{3}\cup\{(0,1,0,0,1,1)\}$&$|\mathcal{O}|=64$&
$\mathcal{C}_{15}^{6}$&$\mathcal{C}_{14}^{2}\cup\{(0,1,1,0,1,0)\}$&$|\mathcal{O}|=128$\\
$\mathcal{C}_{15}^{7}$&$\mathcal{C}_{14}^{2}\cup\{(0,1,0,0,1,1)\}$&$|\mathcal{O}|=48$&
$\mathcal{C}_{15}^{8}$&$\mathcal{C}_{14}^{1}\cup\{(0,1,1,0,0,1)\}$&$|\mathcal{O}|=32$\\
$\mathcal{C}_{15}^{9}$&$\mathcal{C}_{14}^{19}\cup\{(1,1,0,0,0,1)\}$&$|\mathcal{O}|=48$&
$\mathcal{C}_{15}^{10}$&$\mathcal{C}_{14}^{25}\cup\{(0,0,1,1,1,0)\}$&$|\mathcal{O}|=64$\\
$\mathcal{C}_{15}^{11}$&$\mathcal{C}_{14}^{8}\cup\{(1,1,1,0,0,0)\}$&$|\mathcal{O}|=96$&
$\mathcal{C}_{15}^{12}$&$\mathcal{C}_{14}^{5}\cup\{(1,0,0,1,0,1)\}$&$|\mathcal{O}|=32$\\
$\mathcal{C}_{15}^{13}$&$\mathcal{C}_{14}^{1}\cup\{(1,1,1,0,0,0)\}$&$|\mathcal{O}|=96$&
$\mathcal{C}_{15}^{14}$&$\mathcal{C}_{14}^{7}\cup\{(1,1,1,1,1,0)\}$&$\mathcal{D}_{6}$\\
\hline
\end{tabular}
\end{center}
%\end{table}

\newpage

%\begin{table}[h]
%\caption{Cap $\mathcal{K}$ with $|\mathcal{K}|=16$}
\begin{center}
\begin{tabular}{|c|c|c||c|c|c|}
\hline
$\mathcal{K}$&CORRESPONDS TO&STABILIZER&$\mathcal{K}$&CORRESPONDS TO&STABILIZER\\
\hline
\hline
$\mathcal{C}_{16}^{1}$&$\mathcal{C}_{15}^{3}\cup\{(1,0,1,1,1,1)\}$&$\mathcal{D}_{5}$&
$\mathcal{C}_{16}^{2}$&$\mathcal{C}_{15}^{1}\cup\{(1,1,1,0,0,0)\}$&$\mathcal{D}_{4}$\\
$\mathcal{C}_{16}^{3}$&$\mathcal{C}_{15}^{25}\cup\{(0,0,1,1,1,0)\}$&$|\mathcal{O}|=32$&
$\mathcal{C}_{16}^{4}$&$\mathcal{C}_{15}^{14}\cup\{(0,1,1,0,0,1)\}$&$\mathcal{C}_{2}\times \mathcal{C}_{2}\times \mathcal{C}_{2}\times \mathcal{C}_{2}$\\
$\mathcal{C}_{16}^{5}$&$\mathcal{C}_{15}^{12}\cup\{(0,1,0,1,1,1)\}$&$|\mathcal{O}|=128$&
$\mathcal{C}_{16}^{6}$&$\mathcal{C}_{15}^{5}\cup\{(1,1,0,0,0,1)\}$&$|\mathcal{O}|=768$\\
$\mathcal{C}_{16}^{7}$&$\mathcal{C}_{15}^{9}\cup\{(1,1,0,1,1,1)\}$&$|\mathcal{O}|=576$&
$\mathcal{C}_{16}^{8}$&$\mathcal{C}_{15}^{17}\cup\{(1,1,0,0,0,1)\}$&$|\mathcal{O}|=384$\\
$\mathcal{C}_{16}^{9}$&$\mathcal{C}_{15}^{6}\cup\{(0,1,0,0,1,1)\}$&$|\mathcal{O}|=192$&
$\mathcal{C}_{16}^{10}$&$\mathcal{C}_{15}^{21}\cup\{(0,0,1,1,1,0)\}$&$|\mathcal{O}|=144$\\
$\mathcal{C}_{16}^{11}$&$\mathcal{C}_{15}^{7}\cup\{(0,0,1,0,1,1)\}$&$|\mathcal{O}|=32$&
$\mathcal{C}_{16}^{12}$&$\mathcal{C}_{15}^{22}\cup\{(0,1,1,1,1,1)\}$&$|\mathcal{O}|=36$\\
$\mathcal{C}_{16}^{13}$&$\mathcal{C}_{15}^{15}\cup\{(0,1,1,1,1,1)\}$&$|\mathcal{O}|=192$&
$\mathcal{C}_{16}^{14}$&$\mathcal{C}_{15}^{8}\cup\{(1,1,1,0,0,0)\}$&$|\mathcal{O}|=48$\\
\hline
\end{tabular}
\end{center}
%\end{table}

\newpage

%\begin{table}[h]
%\caption{Cap $\mathcal{K}$ with $|\mathcal{K}|=17$}
\begin{center}
\begin{tabular}{|c|c|c||c|c|c|}
\hline
$\mathcal{K}$&CORRESPONDS TO&STABILIZER&$\mathcal{K}$&CORRESPONDS TO&STABILIZER\\
\hline
\hline
$\overline{\mathcal{C}}_{17}^{1}$&$\mathcal{C}_{16}^{7}\cup\{(0,1,1,1,0,1)\}$&$|\mathcal{O}|=576$&
$\overline{\mathcal{C}}_{17}^{2}$&$\mathcal{C}_{16}^{14}\cup\{(0,1,1,0,1,0)\}$&$|\mathcal{O}|=384$\\
$\overline{\mathcal{C}}_{17}^{3}$&$\mathcal{C}_{16}^{14}\cup\{(0,1,0,1,0,1)\}$&$|\mathcal{O}|=720$&
$\overline{\mathcal{C}}_{17}^{4}$&$\mathcal{C}_{16}^{21}\cup\{(0,0,1,1,1,1)\}$&$|\mathcal{O}|=11520$\\
$\overline{\mathcal{C}}_{17}^{5}$&$\mathcal{C}_{16}^{22}\cup\{(1,1,0,1,0,0)\}$&$|\mathcal{O}|=40320$&
$\mathcal{C}_{17}^{6}$&$\mathcal{C}_{16}^{2}\cup\{(0,0,1,0,1,1)\}$&$|\mathcal{O}|=48$\\
$\mathcal{C}_{17}^{7}$&$\mathcal{C}_{16}^{8}\cup\{(1,1,1,0,1,1)\}$&$|\mathcal{O}|=96$&
$\mathcal{C}_{17}^{8}$&$\mathcal{C}_{16}^{1}\cup\{(1,0,1,1,0,0)\}$&$\mathcal{D}_{5}$\\
$\mathcal{C}_{17}^{9}$&$\mathcal{C}_{16}^{2}\cup\{(1,1,0,0,0,1)\}$&$|\mathcal{O}|=64$&
$\mathcal{C}_{17}^{10}$&$\mathcal{C}_{16}^{17}\cup\{(1,1,0,0,0,1)\}$&$\mathcal{D}_{4}$\\
$\mathcal{C}_{17}^{11}$&$\mathcal{C}_{16}^{1}\cup\{(1,0,1,0,0,1)\}$&$\mathcal{C}_{2}\times \mathcal{C}_{2}$&
$\mathcal{C}_{17}^{12}$&$\mathcal{C}_{16}^{7}\cup\{(0,1,1,0,1,0)\}$&$|\mathcal{O}|=2304$\\
$\mathcal{C}_{17}^{13}$&$\mathcal{C}_{16}^{9}\cup\{(0,0,1,1,1,0)\}$&$|\mathcal{O}|=48$&
$\mathcal{C}_{17}^{14}$&$\mathcal{C}_{16}^{22}\cup\{(0,0,1,0,1,1)\}$&$|\mathcal{O}|=2688$\\
\hline
\end{tabular}
\end{center}
%\end{table}
%\begin{table}[h]
%\caption{Cap $\mathcal{K}$ with $|\mathcal{K}|=18$}
\begin{center}
\begin{tabular}{|c|c|c||c|c|c|}
\hline
$\mathcal{K}$&CORRESPONDS TO&STABILIZER&$\mathcal{K}$&CORRESPONDS TO&STABILIZER\\
\hline
\hline
$\overline{\mathcal{C}}^{1}$&$\mathcal{C}_{17}^{14}\cup\{(1,1,1,1,0,0)\}$&$|\mathcal{O}|=10752$&
$\mathcal{C}_{18}^{2}$&$\mathcal{C}_{17}^{21}\cup\{(1,0,1,1,0,0)\}$&$|\mathcal{O}|=128$\\
$\mathcal{C}_{18}^{3}$&$\mathcal{C}_{17}^{22}\cup\{(0,1,1,1,0,0)\}$&$|\mathcal{O}|=64$&
$\mathcal{C}_{18}^{4}$&$\mathcal{C}_{17}^{35}\cup\{(1,0,1,1,0,0)\}$&$|\mathcal{O}|=192$\\
$\mathcal{C}_{18}^{5}$&$\mathcal{C}_{17}^{12}\cup\{(0,1,1,0,0,1)\}$&$|\mathcal{O}|=384$&
$\mathcal{C}_{18}^{6}$&$\mathcal{C}_{17}^{13}\cup\{(0,1,0,1,1,0)\}$&$|\mathcal{O}|=48$\\
$\mathcal{C}_{18}^{7}$&$\mathcal{C}_{17}^{35}\cup\{(1,0,1,0,1,0)\}$&$|\mathcal{O}|=48$&
$\mathcal{C}_{18}^{8}$&$\mathcal{C}_{17}^{8}\cup\{(0,0,1,0,1,1)\}$&$\mathcal{C}_{2}\times \mathcal{C}_{2}$\\
$\mathcal{C}_{18}^{9}$&$\mathcal{C}_{17}^{13}\cup\{(1,1,0,1,1,1)\}$&$|\mathcal{O}|=144$&
$\mathcal{C}_{18}^{10}$&$\mathcal{C}_{17}^{8}\cup\{(0,1,1,0,0,1)\}$&$\mathcal{D}_{4}$\\
$\mathcal{C}_{18}^{11}$&$\mathcal{C}_{17}^{8}\cup\{(1,0,1,0,0,1)\}$&$\mathcal{C}_{2}\times \mathcal{C}_{2}\times \mathcal{C}_{2}$&
$\mathcal{C}_{18}^{12}$&$\mathcal{C}_{17}^{13}\cup\{(1,0,1,1,0,0)\}$&$\mathcal{D}_{4}\times \mathcal{C}_{2}$\\
$\mathcal{C}_{18}^{13}$&$\mathcal{C}_{17}^{9}\cup\{(1,0,0,1,1,0)\}$&$|\mathcal{O}|=64$&
$\mathcal{C}_{18}^{14}$&$\mathcal{C}_{17}^{6}\cup\{(0,1,1,1,0,0)\}$&$\mathcal{D}_{6}$\\
\hline
\end{tabular}
\end{center}
%\end{table}

\newpage

%\begin{table}[h]
%\caption{Cap $\mathcal{K}$ with $|\mathcal{K}|=19$}
\begin{center}
\begin{tabular}{|c|c|c||c|c|c|}
\hline
$\mathcal{K}$&CORRESPONDS TO&STABILIZER&$\mathcal{K}$&CORRESPONDS TO&STABILIZER\\
\hline
\hline
$\mathcal{C}_{19}^{1}$&$\mathcal{C}_{18}^{2}\cup\{(0,1,1,0,1,0)\}$&$|\mathcal{O}|=1152$&
$\mathcal{C}_{19}^{2}$&$\mathcal{C}_{18}^{9}\cup\{(0,1,1,1,0,0)\}$&$|\mathcal{O}|=36$\\
$\mathcal{C}_{19}^{3}$&$\mathcal{C}_{18}^{7}\cup\{(0,0,1,0,1,1)\}$&$\mathcal{S}_{3}$&
$\mathcal{C}_{19}^{4}$&$\mathcal{C}_{18}^{4}\cup\{(0,0,1,0,1,1)\}$&$|\mathcal{O}|=48$\\
$\mathcal{C}_{19}^{5}$&$\mathcal{C}_{18}^{3}\cup\{(0,1,0,1,1,0)\}$&$|\mathcal{O}|=32$&
$\mathcal{C}_{19}^{6}$&$\mathcal{C}_{18}^{15}\cup\{(1,1,1,1,0,0)\}$&$|\mathcal{O}|=9216$\\
$\mathcal{C}_{19}^{7}$&$\mathcal{C}_{18}^{5}\cup\{(0,1,1,1,0,0)\}$&$|\mathcal{O}|=144$&
$\mathcal{C}_{19}^{8}$&$\mathcal{C}_{18}^{3}\cup\{(1,1,1,1,1,0)\}$&$|\mathcal{O}|=192$\\
$\mathcal{C}_{19}^{9}$&$\mathcal{C}_{18}^{17}\cup\{(1,1,0,1,0,0)\}$&$|\mathcal{O}|=192$&
$\mathcal{C}_{19}^{10}$&$\mathcal{C}_{18}^{14}\cup\{(1,0,1,0,1,0)\}$&$|\mathcal{O}|=48$\\
$\mathcal{C}_{19}^{11}$&$\mathcal{C}_{18}^{4}\cup\{(1,0,1,1,1,1)\}$&$|\mathcal{O}|=48$&
$\mathcal{C}_{19}^{12}$&$\mathcal{C}_{18}^{13}\cup\{(0,0,1,0,1,1)\}$&$|\mathcal{O}|=96$\\
$\mathcal{C}_{19}^{13}$&$\mathcal{C}_{18}^{12}\cup\{(1,0,1,1,1,1)\}$&$|\mathcal{O}|=48$&
$\mathcal{C}_{19}^{14}$&$\mathcal{C}_{18}^{14}\cup\{(0,1,1,0,1,0)\}$&$\mathcal{D}_{6}\times \mathcal{C}_{2}$\\
\hline
\end{tabular}
\end{center}
%\end{table}
%\begin{table}[h]
%\caption{Cap $\mathcal{K}$ with $|\mathcal{K}|=20$}
\begin{center}
\begin{tabular}{|c|c|c||c|c|c|}
\hline
$\mathcal{K}$&CORRESPONDS TO&STABILIZER&$\mathcal{K}$&CORRESPONDS TO&STABILIZER\\
\hline
\hline
$\overline{\mathcal{C}}_{20}^{1}$&$\mathcal{C}_{19}^{6}\cup\{(0,1,1,1,1,0)\}$&$|\mathcal{O}|=184320$&
$\mathcal{C}_{20}^{2}$&$\mathcal{C}_{19}^{11}\cup\{(0,1,0,0,1,1)\}$&$|\mathcal{O}|=128$\\
$\mathcal{C}_{20}^{3}$&$\mathcal{C}_{19}^{8}\cup\{(0,1,1,0,1,0)\}$&$|\mathcal{O}|=1536$&
$\mathcal{C}_{20}^{4}$&$\mathcal{C}_{19}^{3}\cup\{(0,0,1,1,0,1)\}$&$|\mathcal{O}|=120$\\
$\mathcal{C}_{20}^{5}$&$\mathcal{C}_{19}^{15}\cup\{(1,0,1,1,1,1)\}$&$|\mathcal{O}|=384$&
$\mathcal{C}_{20}^{6}$&$\mathcal{C}_{19}^{3}\cup\{(1,0,0,1,1,0)\}$&$\mathcal{D}_{4}$\\
$\mathcal{C}_{20}^{7}$&$\mathcal{C}_{19}^{10}\cup\{(0,0,1,1,1,0)\}$&$\mathcal{D}_{6}$&
$\mathcal{C}_{20}^{8}$&$\mathcal{C}_{19}^{8}\cup\{(0,0,1,0,1,1)\}$&$|\mathcal{O}|=64$\\
$\mathcal{C}_{20}^{9}$&$\mathcal{C}_{19}^{2}\cup\{(1,1,0,0,0,1)\}$&$|\mathcal{O}|=144$&
$\mathcal{C}_{20}^{10}$&$\mathcal{C}_{19}^{1}\cup\{(1,1,0,1,1,1)\}$&$|\mathcal{O}|=23040$\\
$\mathcal{C}_{20}^{11}$&$\mathcal{C}_{19}^{3}\cup\{(1,0,1,1,0,0)\}$&$\mathcal{D}_{4}$&
$\mathcal{C}_{20}^{12}$&$\mathcal{C}_{19}^{5}\cup\{(1,1,1,1,1,0)\}$&$|\mathcal{O}|=96$\\
$\mathcal{C}_{20}^{13}$&$\mathcal{C}_{19}^{7}\cup\{(1,1,1,0,0,0)\}$&$|\mathcal{O}|=32$&
$\mathcal{C}_{20}^{14}$&$\mathcal{C}_{19}^{13}\cup\{(1,1,0,1,0,0)\}$&$|\mathcal{O}|=96$\\
\hline
\end{tabular}
\end{center}
%\end{table}

\newpage

%\begin{table}[h]
%\caption{Cap $\mathcal{K}$ with $|\mathcal{K}|=21$}
\begin{center}
\begin{tabular}{|c|c|c||c|c|c|}
\hline
$\mathcal{K}$&CORRESPONDS TO&STABILIZER&$\mathcal{K}$&CORRESPONDS TO&STABILIZER\\
\hline
\hline
$\mathcal{C}_{21}^{1}$&$\mathcal{C}_{20}^{4}\cup\{(1,0,1,1,0,0)\}$&$\mathcal{D}_{5}$&
$\mathcal{C}_{21}^{2}$&$\mathcal{C}_{20}^{3}\cup\{(0,0,1,0,1,1)\}$&$|\mathcal{O}|=192$\\
$\mathcal{C}_{21}^{3}$&$\mathcal{C}_{20}^{6}\cup\{(0,1,1,0,1,0)\}$&$|\mathcal{O}|=32$&
$\mathcal{C}_{21}^{4}$&$\mathcal{C}_{20}^{6}\cup\{(1,0,1,1,0,0)\}$&$\mathcal{D}_{4}\times \mathcal{C}_{2}$\\
$\mathcal{C}_{21}^{5}$&$\mathcal{C}_{20}^{7}\cup\{(1,0,1,1,0,0)\}$&$\mathcal{D}_{6}\times \mathcal{C}_{2}$&
$\mathcal{C}_{21}^{6}$&$\mathcal{C}_{20}^{6}\cup\{(0,1,1,1,1,1)\}$&$|\mathcal{O}|=48$\\
$\mathcal{C}_{21}^{7}$&$\mathcal{C}_{20}^{10}\cup\{(0,1,1,1,0,0)\}$&$|\mathcal{O}|=1920$&
$\mathcal{C}_{21}^{8}$&$\mathcal{C}_{20}^{2}\cup\{(0,1,1,1,1,1)\}$&$|\mathcal{O}|=32$\\
$\mathcal{C}_{21}^{9}$&$\mathcal{C}_{20}^{17}\cup\{(1,1,0,1,1,1)\}$&$|\mathcal{O}|=3072$&
$\mathcal{C}_{21}^{10}$&$\mathcal{C}_{20}^{7}\cup\{(1,1,0,0,0,1)\}$&$|\mathcal{O}|=48$\\
$\mathcal{C}_{21}^{11}$&$\mathcal{C}_{20}^{7}\cup\{(1,0,0,1,1,0)\}$&$\mathcal{D}_{4}\times \mathcal{C}_{2}$&
$\mathcal{C}_{21}^{12}$&$\mathcal{C}_{20}^{15}\cup\{(1,1,0,1,1,1)\}$&$|\mathcal{O}|=1920$\\
$\mathcal{C}_{21}^{13}$&$\mathcal{C}_{20}^{15}\cup\{(0,1,1,1,0,0)\}$&$|\mathcal{O}|=720$&
$\mathcal{C}_{21}^{14}$&$\mathcal{C}_{20}^{12}\cup\{(0,1,1,1,1,1)\}$&$|\mathcal{O}|=1008$\\
$\mathcal{C}_{21}^{15}$&$\mathcal{C}_{20}^{2}\cup\{(0,0,1,0,1,1)\}$&$|\mathcal{O}|=48$&
$\mathcal{C}_{21}^{16}$&$\mathcal{C}_{20}^{3}\cup\{(0,1,0,1,1,0)\}$&$|\mathcal{O}|=384$\\
\hline
\end{tabular}
\end{center}
%\end{table}
%\begin{table}[h]
%\caption{Cap $\mathcal{K}$ with $|\mathcal{K}|=22$}
\begin{center}
\begin{tabular}{|c|c|c||c|c|c|}
\hline
$\mathcal{K}$&CORRESPONDS TO&STABILIZER&$\mathcal{K}$&CORRESPONDS TO&STABILIZER\\
\hline
\hline
$\mathcal{C}_{22}^{1}$&$\mathcal{C}_{21}^{6}\cup\{(0,1,0,0,1,1)\}$&$|\mathcal{O}|=384$&
$\mathcal{C}_{22}^{2}$&$\mathcal{C}_{21}^{3}\cup\{(1,1,0,0,0,1)\}$&$|\mathcal{O}|=192$\\
$\mathcal{C}_{22}^{3}$&$\mathcal{C}_{21}^{2}\cup\{(0,1,0,1,1,0)\}$&$|\mathcal{O}|=48$&
$\mathcal{C}_{22}^{4}$&$\mathcal{C}_{21}^{4}\cup\{(0,1,1,1,1,1)\}$&$|\mathcal{O}|=144$\\
$\mathcal{C}_{22}^{5}$&$\mathcal{C}_{21}^{9}\cup\{(0,1,1,0,0,1)\}$&$|\mathcal{O}|=768$&
$\mathcal{C}_{22}^{6}$&$\mathcal{C}_{21}^{9}\cup\{(1,1,1,1,1,0)\}$&$|\mathcal{O}|=6144$\\
$\mathcal{C}_{22}^{7}$&$\mathcal{C}_{21}^{1}\cup\{(1,0,1,1,1,1)\}$&$\mathcal{D}_{4}\times \mathcal{C}_{2}$&
$\mathcal{C}_{22}^{8}$&$\mathcal{C}_{21}^{1}\cup\{(0,1,0,0,1,1)\}$&$|\mathcal{O}|=120$\\
$\mathcal{C}_{22}^{9}$&$\mathcal{C}_{21}^{14}\cup\{(1,1,1,1,0,1)\}$&$|\mathcal{O}|=2688$&
$\mathcal{C}_{22}^{10}$&$\mathcal{C}_{21}^{2}\cup\{(1,0,1,0,0,1)\}$&$|\mathcal{O}|=3840$\\
$\mathcal{C}_{22}^{11}$&$\mathcal{C}_{21}^{12}\cup\{(0,1,1,1,0,0)\}$&$|\mathcal{O}|=11520$&
$\mathcal{C}_{22}^{12}$&$\mathcal{C}_{21}^{1}\cup\{(1,0,0,1,1,0)\}$&$\mathcal{D}_{4}\times \mathcal{C}_{2}$\\
$\mathcal{C}_{22}^{13}$&$\mathcal{C}_{21}^{5}\cup\{(1,1,0,0,0,1)\}$&$|\mathcal{O}|=72$&
$\mathcal{C}_{22}^{14}$&$\mathcal{C}_{21}^{2}\cup\{(1,0,1,1,1,1)\}$&$|\mathcal{O}|=64$\\
$\mathcal{C}_{22}^{15}$&$\mathcal{C}_{21}^{6}\cup\{(1,1,0,1,0,0)\}$&$|\mathcal{O}|=1008$&&&\\
\hline
\end{tabular}
\end{center}
%\end{table}

%\newpage

%\begin{table}[h]
%\caption{Cap $\mathcal{K}$ with $|\mathcal{K}|=23$}
\begin{center}
\begin{tabular}{|c|c|c||c|c|c|}
\hline
$\mathcal{K}$&CORRESPONDS TO&STABILIZER&$\mathcal{K}$&CORRESPONDS TO&STABILIZER\\
\hline
\hline
$\mathcal{C}_{23}^{1}$&$\mathcal{C}_{22}^{6}\cup\{(1,1,1,1,0,1)\}$&$|\mathcal{O}|=21504$&
$\mathcal{C}_{23}^{2}$&$\mathcal{C}_{22}^{3}\cup\{(1,1,0,1,0,0)\}$&$|\mathcal{O}|=768$\\
$\mathcal{C}_{23}^{3}$&$\mathcal{C}_{22}^{2}\cup\{(1,0,1,1,0,0)\}$&$|\mathcal{O}|=32$&
$\mathcal{C}_{23}^{4}$&$\mathcal{C}_{22}^{2}\cup\{(1,0,1,1,1,1)\}$&$|\mathcal{O}|=48$\\
$\mathcal{C}_{23}^{5}$&$\mathcal{C}_{22}^{5}\cup\{(0,1,1,1,0,0)\}$&$|\mathcal{O}|=1152$&
$\mathcal{C}_{23}^{6}$&$\mathcal{C}_{22}^{1}\cup\{(1,1,1,1,0,1)\}$&$|\mathcal{O}|=384$\\
$\mathcal{C}_{23}^{7}$&$\mathcal{C}_{22}^{4}\cup\{(1,0,1,1,1,1)\}$&$|\mathcal{O}|=48$&
$\mathcal{C}_{23}^{8}$&$\mathcal{C}_{22}^{1}\cup\{(1,0,1,1,0,0)\}$&$|\mathcal{O}|=144$\\
$\mathcal{C}_{23}^{9}$&$\mathcal{C}_{22}^{3}\cup\{(0,1,1,1,1,1)\}$&$|\mathcal{O}|=336$&&&\\
\hline
\end{tabular}
\end{center}
%\end{table}
%\newpage
%\begin{table}[h]
%\caption{Cap $\mathcal{K}$ with $|\mathcal{K}| = 24$}
\begin{center}
\begin{tabular}{|c|c|c||c|c|c|}
\hline
$\mathcal{K}$&CORRESPONDS TO&STABILIZER&$\mathcal{K}$&CORRESPONDS TO&STABILIZER\\
\hline
\hline
$\mathcal{C}_{24}^{1}$&$\mathcal{C}_{23}^{1}\cup\{(0,0,1,1,0,1)\}$&$|\mathcal{O}|=516096$&
$\mathcal{C}_{24}^{2}$&$\mathcal{C}_{23}^{7}\cup\{(1,1,0,1,1,1)\}$&$|\mathcal{O}|=1152$\\
$\mathcal{C}_{24}^{3}$&$\mathcal{C}_{23}^{2}\cup\{(1,0,1,1,1,1)\}$&$|\mathcal{O}|=512$&
$\mathcal{C}_{24}^{4}$&$\mathcal{C}_{23}^{4}\cup\{(0,1,0,1,0,1)\}$&$|\mathcal{O}|=192$\\
$\mathcal{C}_{24}^{5}$&$\mathcal{C}_{23}^{2}\cup\{(1,0,1,0,0,1)\}$&$|\mathcal{O}|=6144$&
$\mathcal{C}_{24}^{6}$&$\mathcal{C}_{23}^{3}\cup\{(1,0,1,1,1,1)\}$&$|\mathcal{O}|=72$\\
$\mathcal{C}_{24}^{7}$&$\mathcal{C}_{23}^{3}\cup\{(0,1,0,0,1,1)\}$&$|\mathcal{O}|=96$&
$\mathcal{C}_{24}^{8}$&$\mathcal{C}_{23}^{1}\cup\{(0,1,1,0,0,1)\}$&$|\mathcal{O}|=2688$\\
\hline
\end{tabular}
\end{center}
%\end{table}
%\begin{table}[h]
%\caption{Cap $\mathcal{K}$ with $|\mathcal{K}|  \geq 25$}
\begin{center}
\begin{tabular}{|c|c|c|}
\hline
$\mathcal{K}$&CORRESPONDS TO&STABILIZER\\
\hline
\hline
$\mathcal{C}_{25}^{1}$&$\mathcal{C}_{24}^{4}\cup\{(0,1,0,0,1,1)\}$&$|\mathcal{O}|=720$\\
$\mathcal{C}_{25}^{2}$&$\mathcal{C}_{24}^{2}\cup\{(0,1,0,0,1,1)\}$&$|\mathcal{O}|=144$\\
$\mathcal{C}_{25}^{3}$&$\mathcal{C}_{24}^{3}\cup\{(0,1,0,0,1,1)\}$&$|\mathcal{O}|=1152$\\
$\mathcal{C}_{25}^{4}$&$\mathcal{C}_{24}^{3}\cup\{(1,0,1,0,0,1)\}$&$|\mathcal{O}|=768$\\
$\mathcal{C}_{25}^{5}$&$\mathcal{C}_{24}^{1}\cup\{(0,1,1,0,0,1)\}$&$|\mathcal{O}|=64512$\\
\hline
\hline
$\mathcal{C}_{26}^{1}$&$\mathcal{C}_{25}^{2}\cup\{(1,1,1,1,0,1)\}$&$|\mathcal{O}|=768$\\
$\mathcal{C}_{26}^{2}$&$\mathcal{C}_{25}^{1}\cup\{(1,0,1,1,0,0)\}$&$|\mathcal{O}|=720$\\
$\mathcal{C}_{26}^{3}$&$\mathcal{C}_{25}^{1}\cup\{(1,0,1,0,0,1)\}$&$|\mathcal{O}|=11520$\\
$\mathcal{C}_{26}^{4}$&$\mathcal{C}_{25}^{4}\cup\{(0,0,1,1,0,1)\}$&$|\mathcal{O}|=18432$\\
\hline
\hline
$\mathcal{C}_{27}^{1}$&$\mathcal{C}_{26}^{1}\cup\{(1,0,1,0,0,1)\}$&$|\mathcal{O}|=1920$\\
$\mathcal{C}_{27}^{2}$&$\mathcal{C}_{26}^{1}\cup\{(1,1,0,1,0,0)\}$&$|\mathcal{O}|=9216$\\
\hline
\hline
$\mathcal{C}_{28}^{1}$&$\mathcal{C}_{27}^{2}\cup\{(0,1,1,0,1,0)\}$&$|\mathcal{O}|=258048$\\
$\mathcal{C}_{28}^{2}$&$\mathcal{C}_{27}^{1}\cup\{(0,1,0,1,0,1)\}$&$|\mathcal{O}|=9216$\\
\hline
\hline
$\mathcal{C}_{29}^{1}$&$\mathcal{C}_{28}^{1}\cup\{(1,0,1,0,0,1)\}$&$|\mathcal{O}|=64512$\\
\hline
\hline
$\mathcal{C}_{30}^{1}$&$\mathcal{C}_{29}^{1}\cup\{(0,1,0,1,0,1)\}$&$|\mathcal{O}|=645120$\\
\hline
\hline
$\mathcal{C}_{31}^{1}$&$\mathcal{C}_{30}^{1}\cup\{(1,1,0,0,0,1)\}$&$|\mathcal{O}|=9999360$\\
\hline
\hline
$\overline{\mathcal{C}}_{32}^{1}$&$\mathcal{C}_{31}^{1}\cup\{(0,0,1,1,0,1)\}$&$|\mathcal{O}|=319979520$\\
\hline
\end{tabular}
\end{center}
%\end{table}

\newpage
All the complete caps of sizes $13$, $17$, $18$ and $20$ have projective frame. The $32$- complete cap has no projective frame.

\end{document}